\documentstyle[12pt,twoside]{article}
\newtheorem{theorem}{Theorem}
\newtheorem{proposition}{Proposition}
\newtheorem{conjecture}{Conjecture}
\newtheorem{lemma}{Lemma}

\newtheorem{remark}{Remark}
\newtheorem{corollary}{Corollary}
\newtheorem{definition}{Definition}

\def\R{I\!\!R}
\def\C{I\!\!\!\!C}
\def\N{I\!\!N}

\def\vsni{\vskip 0.2cm}

\def\ui{[0,1]}

\def\A{{\cal A}}
\def\B{{\cal B}}

\def\be{\begin{equation}}
\def\ee{\end{equation}}
\def\qed{\diamondsuit}

\def\z{\zeta}

\def\b{\beta}

\def\P{{\cal P}}

\def\B{\cal{B}}

\def\di{\displaystyle}

\begin{document}

\title{A one-parameter family of analytic Markov maps with an intermittency transition}
\date{}

\author{Manuela Giampieri\thanks{Dipartimento di Matematica della II Universit\`a
di Roma `Tor Vergata', via della Ricerca Scientifica, I-00133 Roma, Italy. e-mail: $<$giampieri@mat.uniroma2.it$>$} and
Stefano Isola\thanks{Dipartimento di Matematica e Informatica dell'Universit\`a
di Camerino and INFM, via Madonna delle Carceri, 62032 Camerino, Italy.
e-mail: $<$stefano.isola@unicam.it$>$.} }
\maketitle
\begin{abstract} 
\noindent
In this paper we introduce and study a one-parameter family of piecewise analytic 
interval maps having the tent map and the Farey map as extrema. Among other things, we construct
a Hilbert space of analytic functions left invariant by the Perron-Frobenius operator of all these maps and
study the transition between
discrete and continuous spectrum when approaching the intermittent situation.
\end{abstract}
\vskip 1cm
\section{Introduction}
Expanding maps of the unit interval have been widely studied in the last decades in that several problems concerning their statistical behaviour can be  
treated by the powerful technique of transfer operators and thermodynamic formalism \cite{C}, \cite{Ba1}, \cite{Ma2}. 
On the other hand,
in recent years an increasing interest has been carried on maps which are expanding everywhere 
but on a marginally unstable fixed point in a neighbourhood of which trajectories are considerably slowed down 
leading to an interplay of chaotic and regular dynamics characteristic of intermittent systems \cite{PM}, \cite{Sch}. Several
approaches have been proposed to extend the above mentioned techniques to this situation, 
in particular to characterize the nature possible phase transitions \cite{P1}, \cite{PS} 
and that of the spectrum of the transfer operator \cite{Rug}, \cite{Is}. In this paper we 
introduce a one-parameter family of piecewise analytic maps smoothly interpolating between the tent map and the Farey map
and use it to investigate the passage between the uniformly expanding situation and the intermittent one in both perspectives:
that of thermodynamics and that of spectral theory. 

\noindent
The paper is organized as follows:
Section \ref{pre} is devoted to introduce the model and derive some of its properties along with those of an induced version of it.
In Section \ref{free} we discuss large deviation properties and show in particular how the free energy
gets non-analytic in the intermittent limit.
Section \ref{tran} deals with the spectral analysis of the transfer, or Perron-Frobenius, operator.
We construct a Hilbert space of analytic functions where this operator gets a particularly
expressive integral representation (which becomes symmetric in the intermittent limit) and study the mechanism
with which a continuous component in the spectrum (intermittent situation) comes out of a purely discrete 
spectrum (expanding situation).
Finally we extend the above construction to study a family of operator-valued power series by which we
obtain a characterization of the analytic properties
 of the 
dynamical zeta functions \cite{Ba2} for both the original
 map and its induced version discussed in Section \ref{pre}.

\section{Preliminaires} \label{pre}
Let $r\in [0,1]$ be a real parameter and consider the family of piecewise real-analytic maps $F_r$
of the interval $[0,1]$ defined as
\be
F_r(x)= \cases{ F_{r,0}(x),  &if $\; 0\leq x\leq 1/2$\, , \cr
                F_{r,1}(x),  &if $\; 1/2 < x\leq 1$\, ,   \cr }
\ee
where
\be
F_{r,0}(x)={(2-r)x\over 1-rx}\quad\hbox{and}\quad F_{r,1}(x)=F_{r,0}(1-x)={(2-r)(1-x)\over 1-r+rx}\cdot
\ee
Some properties of this family are listed below:

\begin{enumerate}

\item For $r=0$ we find the {\sl tent map}
\be
F_0(x)= \cases{ 2x,  &if $\; 0\leq x\leq 1/2$\, , \cr
                2(1-x),  &if $\; 1/2 < x\leq 1$\, ,   \cr }
\ee
whereas for $r=1$ one has the {\sl Farey map}
\be
F_1(x)= \cases{ {x\over 1-x},  &if $\; 0\leq x\leq 1/2$\, , \cr
                {1-x\over x},  &if $\; 1/2 < x\leq 1$\, .   \cr }
\ee
The latter provides an example of an {\sl intermittent} map  in that
the fixed
point at the origin is neutral (see below).

\item The left branch  $F_{r,0}$ satisfies $F_{r,0}(0)=0$, i.e. the origin is always a fixed point, and
is conjugated to the map $T_r$ defined as
\be
T_r(x)= {x-r\over 2-r}\cdot
\ee
More precisely, we have, for all $n\geq 1$,
\be 
F_{r,0}^n(x) = J^{-1}\circ T_r^n \circ J (x)
\ee
with
\be
J(x)=J^{-1}(x)=1/x.
\ee
Notice that for $r=1$ the map $T_r$ becomes the left translation $x\to x-1$.
Moreover, for each $r\in [0,1]$ there is a unique point $x_1$ left fixed by $F_{r,1}$, i.e. $F_{r,1}(x_1)=x_1$ with 
\be\label{fixpt}
x_1 = {\sqrt{9-4r} -(3-2r) \over 2r}\, \cdot
\ee
\item
The derivative is given by
\be\label{derivative}
F_r'(x)= \cases{\;\;\, F_{r,0}'(x)={2-r\over (1-rx)^2},  &if $\; 0\leq x\leq 1/2$\, , \cr
               F_{r,1}'(x)=-{(2-r)\over  (1-r+rx)^2},  &if $\; 1/2 < x\leq 1$\, ,   \cr }
\ee
so that $F_{r,0}'(x)>0$ and $F_{r,1}'(x)<0$ for all $r\in [0,1]$. 
Therefore $F_r$ is increasing on the interval
$[0,1/2]$ and decreasing on $[1/2,1]$. In addition we have
\be\label{derfix1}
\inf_{x\in [0,1]}|F_r'(x)| = F_{r,0}'(0) =- F'_{r,1}(1)= 2-r=: \rho.
\ee
This means that for $r<1$ the map $F_r$ is uniformly expanding, i.e. $|F_r'|\geq \rho >1$,  thus providing an example of 
{\sl analytic Markov map}
\cite{Ma2}.
In particular we have
\be\label{derfix2}
F'_{r,1}(x_1)={-4\rho \over \left(\sqrt{1+4\rho} -1\right)^2}= {-\left(\sqrt{1+4\rho} +1\right)^2 \over 4\rho}\, \cdot
\ee
On the contrary, for $r=1$ one has $|F'(x)|>1$ for $x>0$ but $F'(0)=1$.
\end{enumerate}

\subsection{Inverse branches and renormalisation}

Set $\rho =2-r$ and 
\be
S_r(x) =T_r^{-1}(x) = \rho x +r
\ee
so that
\be
S_0(x)=2x\quad\hbox{and}\quad S_1(x)=x+1.
\ee
The left inverse branch $\Phi_{r,0}$ of $F_r$ is then given as
\be\label{conj}
\Phi_{r,0}(x) = J^{-1}\circ S_r \circ J (x) = {x\over \rho + rx}=
{1\over 2}-{1\over 2}\left({\rho-\rho\, x\over \rho +rx}\right)\, ,
\ee
whereas the right inverse branch $\Phi_{r,1}$ is
\be\label{conj1}
\Phi_{r,1}(x)=1- \Phi_{r,0}(x)= 1-{x\over \rho + rx}={1\over 2}+{1\over 2}\left({\rho-\rho\, x\over \rho +rx}\right)\cdot
\ee
Note that 
$\Phi_{r,0}(x)$ and $\Phi_{r,1}(x)$ are holomorphic in $H_{-{\rho\over r}}$ and bounded in 
$H_{-{\rho\over r}+\epsilon}$ for all $\epsilon >0$, with  
\be
H_{\alpha} := \{x\in \C : {\rm Re}\, x > \alpha\}.
\ee
Eq. (\ref{conj}) allows us to write an explicit expression  for the iterates $\Phi^n_{r,0}$ of $\Phi_{r,0}$.
\begin{lemma} \label{iterates}
\be
\Phi^n_{r,0}(x) =\left(  {\rho^n\over x} + r\, \sum_{k=0}^{n-1}\rho^k\right)^{-1}
\ee
\end{lemma}
{\sl Proof.} Reasoning inductively in $n$ one obtains 
\be
S^n_r(x) = \rho^n x +r\sum_{k=0}^{n-1}\rho^k
\ee
and the claimed result follows upon applying (\ref{conj}). $\qed$
\vsni
\noindent
Note that (\ref{conj}) can be rewritten as
\be\label{abel}
J(\Phi_{r,0}(x)) = S_r\circ J(x) = \rho \, J(x) + r
\ee
which can be viewed as a generalized Abel equation (see \cite{deB}), to which it actually reduces when $r=1$.
This suggests that $\Phi_{r,0}$ satisfies some non-linear fixed point equation.
Indeed, let ${\cal R}_r$ be the renormalisation operator acting as
\be\label{reno}
{\cal R}_r \Psi (x) = \alpha\, \Psi \left(\Psi \left({x\over\beta}\right)\right),\quad\hbox{with}\quad \alpha = 3-r
\quad\hbox{and}\quad \beta ={3-r\over 2-r}\, \cdot
\ee
\begin{proposition}\label{ren}
For all $r\in [0,1]$ the map $\Phi_{r,0}$ satisfies ${\cal R}_r \Phi_{r,0} = \Phi_{r,0}$ with
boundary conditions $\Phi_{r,0}(0)=0$ and $\Phi_{r,0}'(0)=(2-r)^{-1}=\rho^{-1}$.
\end{proposition}
\noindent
{\sl Proof}. We have
\begin{eqnarray}
{\cal R}_r \Phi_{r,0}(x) &=& \alpha \, \Phi_{r,0}\left(\Phi_{r,0}\left( {x\over\beta}\right) \right)=
\alpha \, \Phi_{r,0}\left( {{x\over\beta}\over \rho + {rx\over \beta}}\right)\nonumber \\
&=& \alpha \, \Phi_{r,0}\left( {x\over \alpha + rx}\right) =  
\alpha \, { {x\over \alpha + rx}\over \rho + {x\over \alpha + rx} } \nonumber \\
&=& \alpha \, { x\over \rho(\alpha + rx)+ rx }= { x\over \rho+ rx }= \Phi_{r,0}(x) \cdot \qquad \qquad \qed \nonumber
\end{eqnarray}
\begin{remark}
For $r=1$ the scaling factors $\alpha$ and $\beta$ in (\ref{reno}) are both
equal to $2$ and the function $\Phi_{0,1}$ is but the fixed point of the Feigenbaum renormalisation equation
with intermittency boundary conditions obtained in {\rm \cite{HR}}.
\end{remark}
\subsection{Invariant measure, induced map, characteristic exponent} \label{im}
We let $\P$ denote the Perron-Frobenius, or transfer operator associated to the map $F_r$ (see \cite{Ba1}). 
It acts on a function $f: [0,1] \to \C$ as
\be
\P f (x) = \sum_{y\, :Ê\, F_r(y)=x} {f (y) \over |F_r'(y)|},
\ee
or, more explicitly,
\be\label{expl}
\P f (x) = {\rho \over (\rho+rx)^2}\left[ f \left({x \over \rho+rx}\right)+
f \left(1-{x \over \rho+rx}\right) \right]\, \cdot
\ee
A fundamental property of this operator is that if there is a measurable function $f$ which satisfies
the fixed point equation 
\be\label{fixed}
\P f (x) = f(x)
\ee
then $f$ is the density of an absolutely continuous measure $\nu_r$ on $[0,1]$ which is $F_r$-invariant,
that is $\nu_r (E) =\nu_r(F_r^{-1}E)$ for all measurable $E \subseteq [0,1]$.

\begin{theorem}\label{invmeas}
The function $e_r(x)=K_r/(1-r+rx)$ where $K_r$ is a given positive constant, is a solution of equation (\ref{fixed}) and thus represents
the density of an absolutely continuous $F_r$-invariant measure $d\nu_r (x) =e_r(x) dx$.
\end{theorem}

\noindent
{\sl Proof}. By virtue of (\ref{expl}) we have
\begin{eqnarray} 
\P \, e_r (x) &=&  {K_r\, \rho \over (\rho+rx)^2}\left[{1 \over 1-r+{rx\over \rho+rx}}+
{1 \over 1-{rx\over \rho+rx}} \right] \nonumber \\
&=& {K_r\, \rho \over (\rho+rx)^2}\left[
{\rho+rx \over \rho +rx-2r+r^2-r^2x+rx}+
{\rho+rx \over \rho} \right]
 \nonumber \\
&=& {K_r \over (\rho+rx)}\left[{\rho \over (1-r)\rho +r\rho x}+1 \right]={K_r\over 1-r+rx} \cdot \qquad \qquad \qed \nonumber
\end{eqnarray}
Notice that 
\be \label{whole}
\nu_r([0,1]) = {K_r\over r} \log\left({1\over 1-r}\right)
\ee
and therefore $\nu_r([0,1])\to \infty$ when $r\nearrow 1$. In order to compare the $F_r$-invariant measure
$\nu_r$ with a probability measure $\mu_r$ invariant w.r.t. to an induced map $G_r$ to be defined below,
we shall choose the value of $K_r$ so that
$\nu_r ([1/2,1])=1$. This renders 
\be
K_r = {r \over \log 2 -\log \rho}\cdot
\ee
In particular we have 
\be
\lim_{r\searrow 0}K_r=2 \qquad \hbox{and}\qquad \lim_{r\nearrow 1}K_r={1\over \log 2} \cdot
\ee
Let ${\cal A}_r=\{A_{n}\}_{n\in \N}$ be the countable partition of
$[0,1]$ whose elements are the intervals
$A_{n}=[c_{n},c_{n-1}]$ with
\be
c_{0}=1\quad\hbox{and}\quad c_{n}=\Phi_{r,0}^n(1),\quad n\geq 1.
\ee
As a corollary of Lemma \ref{iterates} we have the explicit expression
\be\label{cienne}
c_{n}={1-r\over \rho^n -r}\, \cdot
\ee
and in particular
\be
\lim_{r\searrow 0}c_{n}=2^{-n} \qquad \hbox{and}\qquad \lim_{r\nearrow 1}c_{n}={1\over n} \cdot
\ee
Set $X=(0,1]\setminus\{c_{n}\}_{n\in \N}$ and let $\tau : X \to \N$ be the 
{\it first passage time} in the interval $A_{1}$, that is
\be
\tau (x)= 1+ \min \{n\geq 0 \;:\; F_r^n(x)\in A_{1}\;\},
\ee 
so that $A_{n}$ is the closure of the set $\{x\in X\, :Ê\, \tau (x)=n\}$. 
On the other hand, the {\it return time} function $\ell : X \to \N $ 
in the interval $A_{1}$ is given by
\be\label{r}
\ell (x)=\min \{n\geq 1 \;: \; F_r^n(x)\in A_{1}\;\}=\tau \circ F_r(x).
\ee
We now prove the following version of Kac's formula: the $\nu_r$-measure of the
whole interval $[0,1]$, that is (\ref{whole}), equals the conditional expectation of the function $\ell$ on the
interval $A_{1}$ (recall that we have set $\nu_r(A_{1})=1$). 
\begin{lemma}\label{K}
\be\label{Kac}
\int_{A_{1}} \ell (x) \nu_r(dx)=\nu_r([0,1]) = { 1 \over \log 2 -\log \rho} \log\left({1\over 1-r}\right)
\ee
\end{lemma}
{\sl Proof.}
Let $B_{n}={\overline {\{x\in A_{1}\, : \, \ell (x)=n\}} }$. Using (\ref{r}) 
we have that $A_{n}=F_r(B_{n})$. Let us show
that $\nu_r(A_{n})=\sum_{k\geq n}\nu_r(B_{k})$.
Indeed, for $n=1$ we have $1=\nu_r(A_{1})=\sum_{k\geq 1}\nu_r(B_{k})$. Moreover, since $\nu_r$
is $F_r$-invariant, $\nu_r(A_{n})=\nu_r(F_r^{-1}A_{n})=\nu_r(A_{n+1})+\nu_r(B_{n+1})$
and the claim follows by induction. Therefore,
$$
\nu_r([0,1]) = \sum_{n\geq 1}\nu_r(A_{n})= \sum_{n\geq 1}n\cdot \nu_r(B_{n})=
\int_{A_{1}} \ell (x) \nu_r(dx), 
$$
and the last identity in (\ref{Kac}) follows from (\ref{whole}). $\qed$
\begin{remark}
As already remarked, the one-parameter family $F_r$ is well suited to study the transition
from a strongly chaotic behaviour, corresponding to the uniformly expanding situation with $r<1$,
to an intermittent behaviour, corresponding to the the tangent bifurcation
point at $r=1$. One interesting item in this study is the divergence type of the
average duration $<\ell>$ of the {\rm laminar regime} as $r\nearrow 1$ (see, e.g., {\rm \cite{Sch}}). In our situation this
is nothing but the expectation of the return time function $\ell$, and by formula (\ref{Kac})
it diverges logarithmically:
\be
<\ell> \; \sim \; {\log{(1/\delta)} \over \log 2}\quad\hbox{as}\quad \delta \equiv 1-r 
\searrow 0.
\ee
\end{remark}
Now, using (\ref{r}) we can express the expected return time $<\ell>$ as
the expected first passage time w.r.t. an absolutely continuous probability 
measure $\mu_r$ obtained by pushing forward $\nu_r$ with $F_{r,1}$, i.e.
\be
<\ell> = \int_0^1 \tau (x) \, \mu_r (dx) 
\ee
where 
\be
\mu_r (E) =   ((F_{r,1})_*\,\nu_r) (E)=
(\nu_r \circ \Phi_{r,1}) (E) .
\ee
Reasoning as in the proof of Lemma \ref{K} one readily verifies that
the converse relation is
\be
\nu_r (E) = \sum_{n\geq 0} (\mu_r \circ \Phi_{r,0}^n)(E).
\ee
In particular we have $\nu_r (A_{n}) = \sum_{l\geq n} \mu_r ( A_{l})$ and 
$\mu_r (A_{n}) = \mu_r (F_{r,1}(B_{n}))= \nu_r (B_{n})$, where $B_{n}$ is as in the 
proof of Lemma \ref{K}.
The measure $\mu_r$ will play the role of reference probability measure 
in our construction.
If we set $h_r(x) = \mu_r (dx) /dx$
then by the foregoing we have
\be\label{densities}
h_r =  |\Phi_{r,1}'|\cdot e_r\circ \Phi_{r,1},\qquad
e_r = \sum_{k=0}^{\infty}(\Phi_{r,0}^{k})^{\prime} \cdot h_r\circ 
\Phi_{r,0}^{k}\, .
\ee
Using the explicit expressions for $e_r$ and $\Phi_{r,1}$
we get
\be
h_r(x) = {K_r\over 2-r+rx},
\ee
which is  monotone non-increasing with $h_r(0)=K_r/\rho$
and $h_r(1)=K_r/2$. 
\vsni
\noindent
{\sc The induced map.} The measure $\mu_r$ is left invariant by a map $G_r$
obtained from $F_r$ by inducing w.r.t. the first passage time $\tau$.
Indeed, by the above we can write
\be
\mu_r (E) = (\nu_r \circ \Phi_{r,1})(E)
=\sum_{n\geq 0} (\mu_r \circ \Phi_{r,0}^n\circ \Phi_{r,1})(E)= \rho (G_r^{-1}E)
\ee
where $G_r :X\to X$ denotes the map:
\be
x\rightarrow G_r(x) = F_r^{\tau (x)}(x),
\ee
which can be extended to all of $\ui$ as
\be\label{ind}
G_r(x) = G_{r,n}(x)=F_r^n(x)=F_{r,1}\circ F_{r,0}^{n-1}(x)\quad\hbox{if}\quad x\in A_{n}^{\circ}
\quad\hbox{for all}\quad n\geq 1,
\ee
$G_r(0)=G_r(1)=1$ and
\be
\lim_{x\nearrow c_{r,n}}G_r(x)=1,\quad \lim_{x\searrow c_{r,n}}G_r(x)=0,\quad n\geq 1.
\ee
The explicit expression for $G_{r,n}$ can be easily obtained from that of $F_r$ and reads
\be\label{explicit}
G_{r,n}(x)= {\rho \over 1-r+rx}\left(1-{x\over c_{n-1}}\right),
\ee
and
\be\label{explicit'}
G'_{r,n}(x)= -{\rho\, (1-r+rc_{n-1}) \over c_{n-1}(1-r+rx)^2},
\ee
with $c_n$ as in (\ref{cienne}) and $r\in \ui$.

\vsni
\noindent
{\bf Example.} 
The induced maps corresponding to the maps $F_0$ and $F_1$ are the map
$G_0(x)$ such that $G_{0,n}(x)=2(1-2^{n-1}x)$ and
 the Gauss map
$G_1(x)= {\di 1\over \di x}\, ({\rm mod}\, 1)$, respectively.
Their invariant densities are $h_0(x)=1$ and $h_1(x)={\di 1\over \di \log 2}\cdot
{\di 1\over \di (1+x)}\,$.
\vskip 0.1cm
\noindent
We now list some properties of $G_r$ which are relevant for our discussion.
\vskip 0.1cm
\noindent
\begin{proposition}

\noindent
\begin{enumerate}
\item {\sl smoothness property}: $G_{r,n}$ is a real analytic diffeomorphism of $A_{n}$ onto $\ui$;
\item {\sl expanding property}: for all $r\in \ui$
$$ 
\inf_{\scriptstyle x\in \ui }|(G^2_r)'(x)|=|(G^2_r)'(1/2)| =4 \; ;
$$
\item {\sl distortion property}:
$$
\sup_{\scriptstyle x,y,z\in \A_n \atop \scriptstyle n\geq 1} 
\left| {G_r''(x)\over G_r'(y)G_r'(z)}\right| = L < \infty.
$$
\end{enumerate}
\end{proposition}
\noindent
{\it Proof.} Statement 1) is an immediate consequence of the
definition. As for 2) notice that for $r<1$ we have 
$\inf_{\scriptstyle x\in \ui }|(G_r)'(x)|=|G_r'(1)|=\rho > 1$, whereas
$\inf_{\scriptstyle x\in \ui }|(G_1)'(x)|=|G_1'(1)|=1$. Therefore, since $G_r'$ is monotone
decreasing we have $\inf_{\scriptstyle x\in \ui }|(G^2_r)'(x)|=|(G^2_r)'(1/2)|=
|G'_r(1/2)\cdot G_r'(1)|={4\over \rho}\cdot \rho =4$.
To show 3), we first observe that the chain rule yields
\be
{G_r''(x)\over (G_r')^2(x) } =
\sum_{k=0}^{\tau (x)-1}{F_r''(F_r^k(x))\over (F_r')^2(F_r^k(x))}\cdot 
{1\over \prod_{j=k+1}^{\tau (x)-1}F_r'(F_r^j(x))}\, \cdot
\ee
On the other hand, one can easily find a positive constant $C_1$ such that
\be
\sup_{x\in \ui}\, {|F_r''(x)|\over |(F_r')^2(x)|}  \leq C_1.
\ee
Moreover by an easy estimate using (\ref{explicit'}) one can find a constant $C_2\geq 1$ so that 
$C_2^{-1}\, G_r'(y)\leq G_r'(x) \leq C_2\, G_r'(y)$ for any choice of
$x,y\in A_n$ and any $n\geq 1$. Hence, by the mean value theorem
$\prod_{j=0}^{n-1}|F_r'(F_r^j(x))| \equiv |G_{r,n}'(x)|\geq C_2^{-1}\, |A_{n}|^{-1}$
whenever $x\in A_n$. 
The assertion now follows putting together the above inequalities. $\qed$

\vsni
\noindent
These properties yield a uniform bound for the buildup of non-linearity
in the induction process. 
\begin{corollary}
Let $x,y \in \ui$ be such that
$G_r^j(x)$ and $G_r^j(y)$ belong to the same atom $A_{k_{j}}$, for $0\leq j \leq n$
and some $n\geq 1$. 
There is a constant $C>0$, independent of $r$, such that
$$
\left|\log {G_r'(x)\over G_r'(y)}\right| \leq C\,  \left({1\over 2}\right)^{n} \, \cdot
$$
\end{corollary}
{\it Proof.} 
\noindent
Taking $x,y\in A_{k_0}$, let $\eta\in A_{k_0}$ be such 
that $|G_r'(\eta)|=|A_{k_0}|^{-1}$. Then using the distortion property listed above we have
\begin{eqnarray}
\left| \, \log {G_r'(x)\over G_r'(y)}\, \right| 
&=&\left|{G_r''(\xi)\over G_r'(\xi) }\right|\cdot |x-y|
\quad\hbox{for some}\quad \xi \in [x,y] \subseteq A_{k_0} \nonumber \\ 
&=& \left|{G_r''(\xi)\over G_r'(\xi)G_r'(\eta) }\right|\cdot {|x-y|\over |A_{k_0}| } 
\leq L\, {|x-y|\over |A_{k_0}| } \, \cdot \nonumber 
\end{eqnarray}
Now, by the expanding property
we can find a constant $C>0$ such that, under the above hypotheses,
$|x-y| \leq C\,L^{-1} |A_{k_0}| \, \beta^n$ with $\beta = 4^{-{1\over 2}}=1/2$.
$\qed$
\vsni
\noindent
Putting together the above and (\cite{Wal}, Theorem 22(3)) we have the following
\begin{proposition}
 The probability measure $\mu_r$ is the unique
absolutely continuous invariant measure for the dynamical system $(\ui , G_r)$. 
Moreover $(G_r, \mu_r)$ is an exact endomorphism.
\end{proposition}
{\sc Characteristic exponents.} 
Finally we show that the measures $\nu_r$ and $\mu_r$ have the same
characteristic exponent. 
Set 
\be 
\chi_{\nu_r}\,= \int_0^1 \log |F_r'(x)| \, \nu_r (dx) \quad\hbox{and}\quad
\chi_{\mu_r} \, = \int_0^1\log |G_r'(x)|\, \mu_r (dx).
\ee
Then we have 
\begin{proposition}\label{entropies} For all $r\in (0,1)$ we have
\begin{eqnarray}\label{exponents}
\chi_{\nu_r} \, =\, \chi_{\mu_r}  = 
{\log(2-r)\, \log\left({1\over 1-r}\right)\over \log\left({2\over 2-r}\right)} 
&-&\log (4-2r)
 \nonumber \\
&-&{1\over \log\left({2\over 2-r}\right)}\left[{\pi^2\over 6} -\log^2 2 -
2\,{\rm Li}_2\left({1\over 2-r}\right) \right]\nonumber 
\end{eqnarray}
which involves the dilogarithm function ${\rm Li}_2(q)=\sum_{k=1}^\infty {q^n\over k^2}$.
In particular,
$$
\lim_{r\searrow 0}\chi_{\nu_r}\,=\,2\log 2,\qquad \lim_{r\nearrow 1}\chi_{\nu_r}\,=\, {\pi^2\over 6\log 2}\cdot
$$
\end{proposition}
{\sl Proof.} To prove the first identity we write, using (\ref{densities}),
\begin{eqnarray}
 &&\int_0^1 \log |F_r'(x)| \, \nu_r (dx) = \int_0^1 \log |F_r'(x)| \, e_r(x)\, dx  \cr
&=&\int_0^1 \log |F_r'(x)| \, \sum_{k=0}^{\infty}h(\Phi_{r,0}^{k}(x))\cdot 
(\Phi_{r,0}^{k})^{\prime}(x) \, dx 
=\sum_{k=0}^{\infty}\int_0^{c_k}\log |F_r'(F_{r,0}^k(x))|\, h_r(x) \, dx \nonumber \\
&=&\sum_{k=1}^{\infty}\int_{A_k}\prod_{j=0}^{k-1}\log |F_r'(F_{r,0}^j(x))|\, h_r(x) \, dx 
=\sum_{k=1}^{\infty}\int_{A_k}\log |G_{r,k}'(x)|\, h_r(x) \, dx \nonumber \\
&=&\int_0^1\log |G_r'(x)|\, h_r(x) \, dx = \int_0^1\log |G_r'(x)|\, \mu_r (dx). \nonumber
\end{eqnarray}
For the second identity we have, using (\ref{derivative}):
\begin{eqnarray}
\chi_{\nu_r} &=& K_r\left[ \int_0^{1\over 2}\left( \log {\rho \over (1-rx)^2}\right) {dx\over 1-r+rx}+
 \int_{1\over 2}^1\left(\log {\rho \over (1-r+rx)^2}\right){dx\over 1-r+rx}\right]\nonumber \\
&=& \log \rho \cdot \nu_r([0,1]) - 2K_r\left[  \int_0^{1\over 2}{\log (1-rx)\over 1-r+rx}\, dx
+ \int_{1\over 2}^1{\log (1-r+rx)\over1-r+rx}\, dx\right] \nonumber
\end{eqnarray}
By (\ref{Kac}) the first term in the r.h.s equals the first term in the r.h.s. of
(\ref{exponents}). Notice that this term has limits $2\log 2$ and $0$ when $r\searrow 0$ and
$r\nearrow 1$, respectively. Furthermore, we have
$$
2K_r\int_{1\over 2}^1{\log (1-r+rx)\over1-r+rx}\, dx = -{K_r\over r} \log^2\left({2-r\over 2}\right)
=\log\left({2-r\over 2}\right),
$$
and 
$$
2K_r  \int_0^{1\over 2}{\log (1-rx)\over 1-r+rx}\, dx 
= 
 {1\over \log\left({2\over 2-r}\right)}
\left[ {\pi^2\over 6}-\log^22 -2\log 2 \log\left({2-r\over 2}\right)-
2\int_{1\over 2-r}^{0}{\log(1-t)\over t}dt \right].
$$
The claimed formula now follows by putting together the above expressions and
noting that
$$
\int_q^{0}{\log(1-t)\over t}dt = {\rm Li}_2(q)\, .\qquad \qed
$$

\section{Free energy and large deviations}\label{free}

For $r\in [0,1)$ we shall consider the $F_r$-invariant probability measure $p_r$ as well as its 
characteristic (or Lyapunov) 
exponent $\lambda_r$ given by
\be
p_r (\,\cdot \,) = {\nu_r(\,\cdot\,) \over \nu_r([0,1])}\quad\hbox{and}\quad
\lambda_r = {\chi_{\nu_r}\over \nu_r([0,1])}
\ee
respectively. Set moreover
\be
u(x) := \log |F_r'(x)|-\lambda_r,
\ee
and
\be
S_n(x) = \sum_{i=0}^{n-1} u(F_r^i(x)) = \log |(F_r^n)'(x)| - n\, \lambda_r.
\ee
For $\beta\in \R$ and $n\geq 1$ we may then define the {\sl partition function} $Z_n(\beta)$ as
\be
Z_n(\beta)=\int_0^1|(F_r^n)'(x)|^\beta\, p_r (dx) = e^{n\beta\lambda_r}\,\int_0^1e^{\beta \, S_n(x)}\, p_r(dx)
\ee
and consider the sequence of functions
\be\label{freenergy}
f_n(\beta)={1\over n} \log Z_n(\beta)  = \lambda_r \,\beta  + {1\over n} \log \int_0^1e^{\beta \, S_n(x)}\, p_r(dx).
\ee
The limit function
\be\label{freen}
f(\beta) = \lim_{n\to \infty} f_n(\beta)
\ee
is called {\sl free energy function} of the characteristic exponent (see \cite{BR}, \cite{C}, \cite{D}).  
Notice that $f(0) =0$. If we set
\be
\langle A \rangle_\beta := \int_0^1 A(x) { e^{\beta \, S_n(x)} \over \int_0^1e^{\beta \, S_n(x)}\, p_r(dx)}\, p_r(dx)
\ee
so that in particular
\be
\langle A \rangle_0 = \int_0^1 A(x) \, p_r(dx),
\ee
then we get
\be
f'_n(\beta) = \lambda_r + {1\over n}\langle S_n \rangle_\beta\quad\hbox{and}\quad f''_n(\beta) =
{1\over n}\left[\, \langle S^2_n \rangle_\beta-\langle S_n \rangle_\beta^2\, \right].
\ee
Now, for each $r<1$ the transformation $F_r$ is a uniformly expanding Markov map and
using standard arguments one sees that the sequence $\{f''_n(\beta)\, : \, n\geq 1\}$
is uniformly bounded on compact sets. This entails that the free energy function $f(\beta)$ is convex
and $C^1$ with $f'(\beta)$ strictly increasing and given by
\be
f'(\beta) = \lambda_r + \lim_{n\to \infty}{1\over n}\langle S_n \rangle_\beta,\qquad \beta \in \R.
\ee
In particular we have
\be
f'(\beta)|_{\beta =0} = \lambda_r \, \cdot
\ee
One finds moreover that the function $\beta \to \beta \, f'(\beta) - f(\beta)$ is decreasing for $\beta <0$
and increasing for $\beta >0$. This relates
$f(\beta)$ to large deviation properties of the sequence of random variables
$S_n(x)$. To see this, notice that $\int_0^1 S_n(x) \, p_r (dx) =0$ for all $n\geq 1$. 
The ergodic theorem then yields
\be
p_r\, \left( \left\{ x\in \ui \, :Ê\, \lim_{n\to \infty} {S_n(x)\over n}=0\right\}\right) =1
\ee
and therefore for each fixed $\alpha >0$ we have
\be
\lim_{n\to \infty} p_r\, \left( \left\{ x\in \ui \, :Ê\, S_n(x) \geq n\alpha\right\}\right)=0.
\ee
For each $r<1$ the dynamical system $([0,1],F_r,p_r)$ satisfies the assumptions of the large deviation theorem
which says that
(see \cite{C}, \cite{D})
\be\label{ld1}
\lim_{n\to \infty} {1\over n} \log p_r\, \left( \left\{ S_n \geq nf'(\beta)-n\lambda_r\right\}\right) = 
\beta \, f'(\beta) - f(\beta)\qquad (\beta \geq 0)
\ee
and 
\be\label{ld2}
\lim_{n\to \infty} {1\over n} \log p_r\, \left( \left\{ S_n \leq nf'(\beta)-n\lambda_r\right\}\right) = 
\beta \, f'(\beta) - f(\beta)\qquad (\beta \leq 0).
\ee
In other words, the probability of finding a deviation of $S_n/n$ 
from its average value $0$ decays exponentially with $n$. Notice that if we set $\alpha =f'(\beta)-\lambda_r$
then the r.h.s. in (\ref{ld1}) and (\ref{ld2}) can be viewed as a
Legendre transform
\be
\phi (\alpha ) =  \beta\, \alpha - (f(\beta) - \beta \lambda_r),\qquad \alpha =f'(\beta)-\lambda_r.
\ee
We now derive upper and lower bounds for the free energy function. First, from  $|F_r'|\geq \rho >1$ it follows
immediately that for $\beta \leq 0$ we have $f(\beta) \leq \beta \log \rho$, the opposite inequality being valid
for positive values of $\beta$. In addition, using either the convexity of $f(\beta)$ or directly (\ref{freenergy})
we get the inequality $f(\beta) \geq \lambda_r\, \beta$, which is valid for all $\beta \in \R$.
Whence, expressing $\lambda_r$ by means of  Lemma \ref{K} and Proposition \ref{entropies} we get
\begin{lemma}\label{bounds} For all $\beta \leq 0$ and $r\in [0,1)$ we have 
$$
\beta\log \rho \geq {f(\beta)} \geq \beta\log \rho + \beta \, \gamma_r
$$
where $\gamma_r\geq 0$ is given by
$$
\gamma_r=
{1\over \log (1-r)}\left[{\pi^2\over 6}-
2\,{\rm Li}_2\left({1\over 2-r} \right)+\log \left({2\over 2-r}\right) \log (4-2r) -\log^2 2 \right]\, .
$$
\end{lemma}
Notice that $\lim_{r\searrow 0}\log \rho=\log 2$ whereas $\lim_{r\nearrow 1}\log \rho=0$, moreover
$\lim_{r\searrow 0}\gamma_r=0$ and $\lim_{r\nearrow 1}\gamma_r=0$. 
\noindent
Therefore from Lemma \ref{bounds} we obtain the
\begin{corollary}
$$
\lim_{r\searrow 0}f(\beta)=\beta \, \log 2,\qquad \lim_{r\nearrow 1}f(\beta) =0 , \qquad (\beta \leq 0).
$$
\end{corollary}
\begin{remark}
This result shows that for $r=1$ the free energy has a discontinuity in its first derivative at $\beta =0$.
This can be interpreted in thermodynamic language as a second order phase-transition \cite{P1}, \cite{FKO}.
\end{remark}
\section{Transfer operators}\label{tran}
The transfer operator $\P$ associated to 
the map $F_r$ has already been introduced in (\ref{expl}) and we write 
\be
\P f (x) = ({\P}_0+{\P}_1)f (x)
\ee
with (recall that $\rho \equiv 2-r$)
\be
{\P}_0f (x)={\rho \over (\rho+rx)^2}\cdot f \left({x \over \rho+rx}\right)\quad\hbox{and}\quad
{\P}_1f (x) ={\rho \over (\rho+rx)^2}\cdot f \left(1-{x \over \rho+rx}\right) \, .
\ee
Several interesting properties of the dynamics generated by $F_r$ are intimately related to ${\rm sp}\,(\P)$,
the spectrum of $\P$ (see \cite{Ba1}). However, the latter depends crucially on the function space 
$\P$ is acting on, which is in general a Banach space. For 
smooth uniformly expanding maps and Banach spaces of sufficiently
regular functions, e.g. the space ${\cal C}^k$ of $k$-times differentiable functions on $\ui$ with $k\geq 0$, the transfer operator is
{\sl quasi-compact}. This means that ${\rm sp}\, (\P)$ is made out of a finite or at most countable set of isolated eigenvalues
with finite multiplicity (the discrete spectrum) and its complementary, the essential spectrum. 
It has been proved in \cite{CI} that for piecewise ${\cal C}^\infty$ expanding Markov maps of the unit interval the essential spectrum 
of $\P$ when acting on ${\cal C}^k$ is a disk of radius 
\be
r_{\rm ess}(\P)\, =\, \exp{f(-k)}\, ,
\ee where $f(\b)$ is the free energy function
discussed in the previous Section. Putting together this result and Lemma \ref{bounds} we obtain the following

\begin{theorem}\label{esspect} For $r\in [0,1)$ the essential spectrum of $\P:{\cal C}^k\to {\cal C}^k$ with $k\geq 0$ is a disk
of radius
$$
 e^{-k\, (\log \rho + \gamma_r)}  \;  \leq \; r_{\rm ess}(\P) \; \leq \;e^{-k\, \log \rho}\, .
$$
\end{theorem}
The above bounds along with standard semicontinuity arguments \cite{K} yield the following
\begin{corollary}\label{one} 
For $r=1$ and for each fixed $k\geq 0$ the essential spectrum of $\P:{\cal C}^k\to {\cal C}^k$
is the unit disk.
\end{corollary}
\subsection{An invariant Hilbert space}
From the above discussion it follows that if we want to understand the nature of the spectrum lying under the `essential spectrum rug'
we have to let $\P$ acting on increasingly smooth test functions as $r$ approaches $1$. In particular, Corollary \ref{one} suggests
that for $r=1$ one should consider suitable spaces of analytic functions. 
In the following definition we shall introduce a Hilbert space of analytic functions which will be shown to be
invariant under $\P$ for each $r\in [0,1]$.

\begin{definition} We denote by ${\cal H}$ 
the Hilbert space of all complex-valued functions $f$ which can be represented as a generalized 
Borel transform
\be\label{represent2}
f (x)= ({\B}\, [\varphi])(x):={1\over x^2}\int_0^\infty  e^{-{t\over x}}\,e^t
\varphi (t)\, dm(t),\quad 
\varphi \in L^2(m),
\ee
with inner product
\be
(f_1,f_2) = \int_0^\infty  \varphi_1(t)\, {\overline {\varphi_2  (t)}}\, dm(t)\quad\hbox{if}\quad f_i={\B}\, [\varphi_i],
\ee
and measure
\be
dm(t) = t\, e^{-t}\, dt.
\ee
\end{definition}
\begin{remark}
An alternative representation can be obtained by a simple change of variable when $x$ is real and positive:
\be\label{borel}
f(x) = {1\over x}\int_0^\infty ds\, e^{-s}\, \psi (sx)\quad \hbox{with}\quad
\psi (t) = t\, \varphi (t)\, .
\ee
Note that a function $f \in {\cal H}$ is analytic in the disk (here $z=x+iy$)
\be
D_1=\{z\in \C : {\rm Re}\, {1\over z} > {1\over 2}\}=\{z\in \C : |z-1|<1\}.
\ee
A Hilbert space identical to ${\cal H}$ apart from a slightly different choice of the
measure $m$ was introduced in \cite{Is} to study the operator $\P$ for $r=1$ (the Farey map), whereas a generalized version of
${\cal H}$ has been used by Prellberg in \cite{P2} to study 
the spectrum of the operator $\P_\beta f (x) = \left({1 \over 1+x}\right)^{2\beta} \left[ f \left({x \over 1+x}\right)+
f \left({1 \over 1+x}\right) \right]$ (thus again for the case $r=1$). 
\end{remark}

\begin{remark}\label{densi}
The invariant densities $e_r(x)=K_r/(\delta+rx)$ and $h_r(x)=K_r/(\rho+rx)$ can be represented as
\be
e_r ={\B}[\phi_r]\quad \hbox{and}\quad h_r ={\B}[\psi_r]
\ee
with
\be\label{densita}
\phi_r (t) = {K_r\over r}\left({1-e^{-{r\over \delta \,}t}\over t}\right) \quad \hbox{and}\quad
\psi_r (t) = {K_r\over r}\left({1-e^{-{r\over \rho \,}t}\over t}\right) ,
\ee
respectively. For the limiting values $r=1$ and $r=0$ we get
\be
\phi_1(t)\,=\,{1\over t\, \log 2} ,\qquad \psi_1(t)\,=\,{1-e^{-t}\over t\,\log 2}\, ,
\ee
and
\be
\phi_0(t)=\lim_{r\searrow 0}\phi_r(t)\,=\,2 , \qquad \psi_0(t)=\lim_{r\searrow 0}\psi_r(t)\,=1\, .
\ee
We point out that $\phi_1$ is not in $L^2(m)$.
\end{remark}
\begin{lemma} \label{pi0}
The space ${\cal H}$ is invariant for ${\P}_0$ and we have
\be
{\P}_0{\B}\, [\varphi] = {\B}\, [M_r\varphi]\, ,
\ee
where $M_r :L^2(m)\to L^2(m)$ is defined as
\be
M_r \varphi (t) = {1\over \rho}\,e^{-{r\over \rho}t}\, \varphi\left({t\over \rho}\right).
\ee 
\end{lemma}
{\sl Proof.} 
\begin{eqnarray}
({\P}_0{\B}\, [\varphi])(x) &=& {\rho \over (\rho +rx)^2}{\B}\, [\varphi]\left({x\over \rho +rx}\right) \nonumber \\
&=&{\rho\over x^2}\int_0^\infty e^{-{\rho +rx\over x}t}\, \varphi (t) \, t\, dt \nonumber \\
&=& {1\over x^2}
\int_0^\infty e^{-{s\over x}}\, {1\over \rho}\,e^{-{r\over \rho}s}\,\varphi 
\left({s\over \rho}\right) \, s\, ds  \nonumber \\ 
&=& {1\over x^2}
\int_0^\infty e^{-{s\over x}}\, e^s \left(M_r\varphi\right) (s)
 \, dm(s) = ({\B}\, [M_r\varphi])(x)\, . \qquad \qquad \qed \nonumber
\end{eqnarray}
The following lemma will instead specify the action of $\P_1$ on ${\cal H}$.
\begin{lemma} \label{pi1} We have
\be
{\P}_1{\B}\, [\varphi] = {\B}\, [ N_r\varphi]\, ,
\ee
where $N_r:L^2(m)\to L^2(m)$ is the operator acting as 
\be
N_r \varphi (t) =  {1\over \rho}\,e^{{\delta\over \rho}t}\, \int_0^\infty 
{J_1\left({2}\sqrt{{st/ \rho}}\right)\over\sqrt{ st/\rho}}\, \varphi (s) \, dm(s)
\ee
where $J_p$ denotes the Bessel function of order $p$.
\end{lemma}
{\sl Proof.} 
\begin{eqnarray}
({\P}_1{\B}\, [\varphi])(x) &=& {\rho \over (\rho +rx)^2}\,{\B}\, [\varphi]\left(1-{x\over \rho +rx}\right)
 \nonumber \\
&=&{\rho\over (\rho -\delta \,x)^2}\int_0^\infty e^{-{\rho +rx\over \rho-\delta x}t} 
\, \varphi (t) \, t\, dt \nonumber \\
&=& {\rho\over (\rho -\delta \,x)^2}
\int_0^\infty e^{-{ \rho -\delta x +\delta x +rx\over \rho -\delta  x}t} \varphi ( t) \,t\, dt  \nonumber \\ 
&=& {\rho\over x^2}{1\over \left({\rho\over x} -\delta\right)^2}
\int_0^\infty e^{-{ t\over {\rho\over x} -\delta }}  e^{- t} \varphi (t)\,t \, dt  \nonumber \\
&=& {\rho \over x^2}  \int_0^\infty e^{-\left({\rho\over x} -\delta\right)t}\left(
\int_0^\infty  J_1(2\sqrt{st})\, \sqrt{t\over s}  \, 
\varphi (s) \,s\,e^{- s}\, ds   \right) \, dt \nonumber \\
&=& {1 \over x^2}  \int_0^\infty e^{-{t\over x}}\, e^{\,{\delta\over \rho}t}
\left(
\int_0^\infty  J_1\left({2}\sqrt{st\over \rho}\right)\, \sqrt{t\over \rho s} \,
\varphi (s) \, dm(s )  \right) \, dt  
\nonumber \\ 
&=& {1 \over x^2}  \int_0^\infty e^{-{t\over x}}\, e^t\, {1\over \rho}\, e^{\,{\delta\over \rho}t}
\left(
\int_0^\infty 
{J_1\left({2}\sqrt{{st/ \rho}}\right)\over\sqrt{s t/\rho}}\, \varphi (s) \, dm(s) \right) \, dm(t)  
\nonumber \\ 
&=& ({\B}\, [N_r\varphi])(x)\, , \nonumber
\end{eqnarray}
where we have used the identity \cite{GR}
\be\label{bessel}
{1\over u^{p+1}}\int_0^\infty  e^{-t/u} \psi (t)\, dt  = \int_0^\infty e^{-tu}\left(
\int_0^\infty \left({t\over s}\right)^{p\over 2}\, J_p(2\sqrt{st})\, \psi (s)\, ds\right) \, dt
\ee
with $p=1$, $u={\rho\over x} -\delta$ and $\psi (s) =s\, e^{- s} 
\varphi (s)$. Finally, by the estimates (see \cite{Er}) $J_1(x) \sim {x/2}$ as $x\searrow 0$ and
$J_1(x) = {\cal O}(x^{-{1\over 2}})$ as $x\to \infty$, along with the inequality $2\delta/\rho <1$, which holds for all $r\in [0,1]$, one easily checks that
$N_r: L^2(m)\to L^2(m)$.
$\qed$
\vsni
\noindent
We can summarize the above in the following
\begin{theorem}\label{pi}
The space ${\cal H}$ is invariant for ${\P}$ and we have
\be\label{Pop}
{\P}{\B}\, [\varphi] = {\B}\, [(M_r+N_r)\varphi]\, .
\ee
\end{theorem}

\subsection{The spectrum of $M_r$ and $\P_0$}
We are now going to study the operator $M_r$. First,
note that it reduces to $M_0\, \varphi(t)=(1/2)\varphi(t/2)$  when $r=0$, whereas for $r=1$ 
yields the multiplication operator $M_1\varphi (t) =e^{-t}\varphi (t)$.
Moreover, for $r\in [0,1)$ its iterates are given by
\be\label{iterate}
M_r^n \varphi (t) = \left(\prod_{k=1}^n e^{-{r\over \rho^k}t}\right)\, \varphi \left({t\over \rho^n}\right)
= {1\over \rho^n}\,e^{-{r\over \delta}t}\, e^{\, {r\over \delta}{t\over \rho^n}}\,  \varphi \left({t\over \rho^n}\right).
\ee
Assume that $\varphi$ is analytic in some open neighbourhood of $t=0$ and
satisfies the eigenvalue equation 
\be\label{autova}
M_r \varphi (t)= {1\over \rho}\,e^{-{r\over \rho}t}\varphi \left({t\over \rho}\right)= \lambda \varphi(t).
\ee
The above equation at $t=0$ writes
$\varphi(0) = \lambda \,\rho\, \varphi(0)$, so that if $\varphi (0)\ne 0$ then $\lambda =1/\rho$. 
In this case, by (\ref{iterate}) we get
\be
\varphi (t) = \rho^n\,M_r^n \varphi (t)\to e^{-{r\over \delta}t}\, \varphi (0)\quad\hbox{as}\quad n\to \infty,
\ee
and therefore
$\varphi (t) = e^{-{r\over \delta}t}\, \varphi (0)$.
If instead $\varphi (0)= 0$ we differentiate (\ref{autova}) to get
\be
-{r\over \rho^2}e^{-{r\over \rho}t}\varphi (t/ \rho)+{e^{-{r\over \rho}t}\varphi' (t/ \rho)\over \rho^2} = \lambda \varphi'(t)
\ee
which at $t=0$ writes $\varphi' (0)/ \rho^2 = \lambda \varphi'(0)$.
Therefore if $\varphi' (0)\ne 0$ then $\lambda =1/\rho^2$.
In this case, differentiating (\ref{iterate})-(\ref{autova})  we get
\be
{1\over \rho^{2n}}\varphi' (t) =  {r\over \delta\,\rho^n}\left({1\over \rho^n}-1\right)e^{-{r\over \delta}t}\, e^{\, {r\over \delta}{t\over \rho^n}}\,  
\varphi \left({t\over \rho^n}\right)+
{e^{-{r\over \delta}t}\, e^{\, {r\over \delta}{t\over \rho^n}}\over \rho^{2n}}  \varphi' \left({t\over \rho^n}\right)
.
\ee
Taking
the limit $n\to \infty$ and noting that $\lim_{n\to \infty} \rho^n \varphi (t/ \rho^n) = t\, \varphi'(0)$ we obtain
\be
\varphi' (t) = \varphi'(0)\, e^{-{r\over \delta}t}\left(1-{r\over \delta}\, t\right),
\ee
which upon integration renders
\be
\varphi (t) = \varphi'(0)\,t\, e^{-{r\over \delta}t}.
\ee
Iterating this argument we have that if $\varphi$ satisfies (\ref{autova}) with $\varphi^{(l)}(0) = 0$ for 
$0\leq l<k-1$ but $\varphi^{(k-1)}(0) \ne 0$ for some $k \geq 1$ then $\lambda =\rho^{-k}$ and 
$\varphi (t) =t^{k-1}\, e^{-{r\over \delta}t}$ (up to a non-zero but otherwise arbitrary constant multiplicative factor).
Denoting by $\Vert \; \Vert_2$ the norm in $L^2(m)$ we also have
\be\label{norm}
\Vert t^{k-1}\, e^{-{r\over \delta}t} \Vert_2^2 = \left({\delta \over 1+r}\right)^{k}(2k-1)!
\ee
It is not hard to see that for all $r\in [0,1)$ the sequence $\{t^{k-1}\, e^{-{r\over \delta}t}\}_{k=1}^\infty$
is a linearly independent family in $L^2(m)$, and by adapting (\cite{He}, p.62, Thm.8) we have that
the linear span of this family is dense in $L^2(m)$. 
Putting together the above along with standard arguments we have proved the following
\begin{proposition} \label{spettro}
For all $r\in [0,1)$ the operator $M_r :L^2(m)\to L^2(m)$ is compact and its spectrum is
given by ${\rm sp}\, (M_r) = \{\mu_k\}_{k\geq1}\cup\{0\}$ with
\be\label{e1}
\mu_k ={1\over \rho^{k}} \equiv \left( \Phi'_{r,0}(0)\right)^k.
\ee
Each eigenvalue $\mu_k$ is simple
and the corresponding (normalized) eigenfunction $\varphi_k$ is given by 
\be\label{eigenfct}
\varphi_k(t)=A_k\, t^{k-1}\, e^{-{r\over \delta}t}
\ee
with
\be
A_k =\left({1+r\over \delta}\right)^{k}{1\over \sqrt{(2k-1)!}}\, \cdot
\ee
\end{proposition}
\begin{remark}\label{tr}
For each fixed $r\in [0,1)$ the operator $M_r :L^2(m)\to L^2(m)$ is actually trace-class. Its trace is easily computed:
\be\label{trM}
{\rm tr}\, M_r = \sum_{k=1}^\infty \rho^{-k} = {1 \over \delta} \, ,
\ee
and satisfies
\be
\lim_{r\searrow 0}{\rm tr}\, M_r\,=\,1\,, \qquad \lim_{r\nearrow 1}{\rm tr}\, M_r\,=\,\infty \, .
\ee
\end{remark}
We recall that for $r=1$ we get the multiplication operator $M_1$ which is self-adjoint in $L^2(m)$, its spectrum is continuous and
given by the closure of the range of the multiplying function, that is the interval $[0,1]$ (see, e.g., \cite{DeV}).
By the above Proposition we see how the continuous spectrum is approached as $r\nearrow 1$: having fixed an
interval $[a,b]\subseteq (0,1]$ an easy computation yields
\be
\# \{ \mu_k \in [a,b]\}\sim {\log \left({b-a\over ab}\right) \over 1-r}
\quad\hbox{as}\quad r\nearrow 1.
\ee
Moreover a simple calculation gives
\be
{\B}\, [t^{k-1}\, e^{-{r\over \delta}t}]= {k! \,\delta^{k+1}\, x^{k-1}\over (\delta + rx)^{k+1}}.
\ee
Therefore by the above and Lemma \ref{pi0} we have the following 
\begin{corollary} \label{P0}
The spectrum of ${\P}_0$ when acting upon ${\cal H}$ is given by 
${\rm sp}\,({\P}_0)={\rm sp}\, (M_r)$. For $r\in [0,1)$ each
eigenvalue $\mu_k$, with $k\geq 1$, is simple
and the corresponding (normalized) eigenfunction ${\chi}_k$ is given by 
\be
{\chi}_k (x)=({\B}_r\, [\varphi_k])(x)=\, {(1+r)^{k} \, k !\, \delta\over \sqrt{(2k -1)!}}\, 
{x^{k-1}\over (\delta + rx)^{k+1}}\, \cdot
\ee
\end{corollary}
\subsection{The spectrum of $N_r$ and $\P_1$}
We now turn to study the action of the operator $N_r$ on $L^2(m)$. We first note that $N_r$ can be viewed as the
composition of the multiplication operator 
\be
\varphi (t) \to {1\over \rho}\, e^{{\delta\over \rho}t}\,\varphi (t),
\ee
(which reduces to the identity for $r=1$) 
and the symmetric integral operator
\be\label{symm}
\varphi (t) \to  \int_0^\infty 
{J_1\left({2}\sqrt{{st/ \rho}}\right)\over\sqrt{{st/ \rho}}}\, \varphi (s) \, dm(s).
\ee
Observing that the (associated) Laguerre polynomials $L_{k}^{1}$ given by \cite{GR}
\be
L_{k}^{1}(s) = {e^s\, s^{-1}\over k!}{d^k\over ds^k} \left(e^{-s}\, s^{k+1}\right)=
\sum_{l=0}^k {k+1 \choose k-l} {(-s)^l\over l!} 
\ee
 form a complete orthogonal basis in $L^2(m)$ and
expanding the kernel of the integral operator defined above on this basis we get 
\be
{J_1\left({2}\sqrt{{st/ \rho}}\right)\over\sqrt{{st/ \rho}}} = 
\sum_{k=0}^\infty L_{k}^{1}(s) \, \, {t^{k}\, e^{-{t\over\rho}} \over \rho^{k}(k+1)!} 
\ee
from which we see that $N_r$ has the representation 
(we keep using the symbol $(\varphi_1,\varphi_2)$ to denote the inner product in $L^2(m)$ as well)
\be
N_r\varphi = \sum_{k=1}^\infty   (\varphi,e_k)\, f_k
\ee
where $e_k, f_k \in L^2(m)$ are given by
\be
e_k(t)={L_{k-1}^{1}(t)}\quad\hbox{and}\quad f_k(t) ={N_r e_k(t)} ={t^{k-1}\, e^{-{r\over\rho}t} \over \rho^{k}\,k!}\,\cdot
\ee
We find
\be
\Vert e_k\Vert_2 =\sqrt{k}\quad\hbox{and}\quad\Vert f_k\Vert_2 = { \sqrt{(2k-1)!}\over  (2+r)^{k}\, k!},
\ee
and therefore, for all $r\in (0,1]$,
\be
\sum Œ  \Vert e_k\Vert_2 \, \Vert f_k\Vert_2 < \infty
\ee
showing that the operator $N_r$ is nuclear in $L^2(m)$. Its trace can be computed as
\be\label{trN}
{\rm tr} \, N_r = \int_0^\infty {e^{-{s\over \rho}}\over \sqrt{\rho}}\, J_1\left( {2s\over \sqrt{\rho}} \right) ds
={1\over 2}\left(1-{1\over\sqrt{1+4\rho}}\right) ,
\ee
with limit values
\be
\lim_{r\searrow 0}{\rm tr}\, N_r\,=\,{1\over 3}\quad\hbox{and}\quad\lim_{r\nearrow 1}{\rm tr}\, N_r\,=\,{\sqrt{5}-1\over 2\sqrt{5}}\, \cdot
\ee
Moreover we find
\be
{\rm tr} \, N^2_r = {1\over \rho} \int_0^\infty \int_0^\infty e^{-{t+s\over \rho}}\left[ J_1\left({2}\sqrt{{st/ \rho}} \right) \right]^2\, 
ds\, dt
={1\over 2}\left({ 1+2\rho\over \sqrt{1+4\rho}}-1\right).
\ee
Note that ${\rm tr} \, N^2_r\leq {\rm tr} \, N_r$, with strict inequality unless $r=0$ where ${\rm tr} \, N_0={\rm tr} \, N^2_0=1/3$.
This suggests that the eigenvalues of $N_r$ are alternately positive and negative. 
Indeed we have the following 
\begin{proposition} \label{spettro2}
For all $r\in [0,1]$ the spectrum of the operator $N_r :L^2(m)\to L^2(m)$ is
given by ${\rm sp}\, (N_r) = \left\{{\nu}_k\right\}_{k= 1}^\infty\cup\{0\}$ with
\be
{\nu}_k= (-1)^{k-1}\left({4\rho\over (1+\sqrt{1+4\rho})^2}\right)^k \equiv 
-\left( \Phi'_{r,1}(x_1)\right)^k.
\ee 
Each eigenvalue is simple
and the corresponding (normalized) eigenfunction ${\psi}_k$ is given by 
\be\label{eigenfct}
{\psi}_k(t)=B_k\, L_{k-1}^1\left(\alpha_r t\right)\, e^{- \beta_r\,t},
\ee
where $\alpha_r={\sqrt{1+4\rho}\over \rho}$, $\beta_r = {1+\sqrt{1+4\rho}\over 2\rho}-1$ and
$$
B_k = {\sqrt{1+4\rho}\over \rho \, \sqrt{k}} \left( 1-{\delta\over \sqrt{1+4\rho}}\right)^k \left[ \sum_{j=0}^{k-1} 
{k\choose j} {k-1 \choose j}
\left({\delta^2\over1+4\rho}\right)^{j} \right]^{-{1\over 2}}\cdot
$$
\end{proposition}
\begin{remark} Note that the eigenvalues of $N_r$ can be written in terms of the function $\beta_r$ as
\be\label{e2}
{\nu}_k={(-1)^{k-1}\over \rho^k} \left({1\over 1+\beta_r}\right)^{2k}\, \cdot
\ee
In particular, $\lim_{r\searrow 0}\beta_r=0$ so that for $r=0$ we have ${\nu}_k=(-1)^{k-1}(2)^{-k}$,
moreover $\lim_{r\nearrow 1}\beta_r={\sqrt{5}-1\over 2}$ and thus for $r=1$ we see that 
${\nu}_k=(-1)^{k-1}\left({\sqrt{5}-1\over 2}\right)^{2k}$.
\end{remark}
\noindent
{\sl Proof of Proposition \ref{spettro2}.} The proof is based on the following Hankel transform (see \cite{Er}, vol 2)
\be\label{hankel}
\int_0^\infty x^{p+{1\over 2}}\, e^{-bx^2}\, L_k^p(ax^2) \, J_p(xy)\, \sqrt{xy}\, dx= 
{(b-a)^k\, y^{p+{1\over 2}}\over 2^{p+1}\, b^{p+k+1}}\, e^{-{y^2\over 4b}}\, L_{k}^p\left[{ay^2\over 4b(a-b)}\right],
\ee
which for $p=1$ can be rewritten in terms of the operator $N_r$ as 
\be\label{formulona}
N_r \left[ e^{-\left({2b\over \sqrt{\rho}} -1\right)t}\, L_{k}^1\left({2at\over \sqrt{\rho}}\right)\right] =
{(b-a)^k \over 4 \, b^{k+2}} \, e^{-\left({1\over 2b\sqrt{\rho}} -{\delta \over \rho}\right)t}\,
  L_{k}^1\left[{at\over 2\sqrt{\rho}\, b(a-b)}\right].
\ee
To make the above identity an eigenvalue equation the following relations have to be satisfied
\be
{2b\over \sqrt{\rho}} -1 = {1\over 2b\sqrt{\rho}} -{\delta \over \rho}\quad\hbox{and}\quad 
{2a\over \sqrt{\rho}} = {a\over 2\sqrt{\rho}\, b(a-b)}.
\ee
This renders
\be
a=b+{1\over 4b}\quad\hbox{and}\quad b={1\pm \sqrt{1+4\rho} \over 4\sqrt{\rho}}.
\ee
It is now easy to check that the only solution giving a function in $L^2(m)$ is that with $b={(1+ \sqrt{1+4\rho}) / 4\sqrt{\rho}}$
and this choice yields the eigenvalues and the eigenfunctions given in the Proposition. Moreover, a standard evaluation of
\be\label{norm2}
\Vert L_{k-1}^1\left(\alpha(r) t\right)\, e^{- \beta(r)\,t} \Vert_2^2 = {\rho^2\over 1+4\rho}
\int_0^\infty \left[ L_{k-1}^1(s)\right]^2\, s\, e^{-\left(1-{\delta\over\sqrt{1+4\rho}}\right)s}\,ds
\ee
yields the claimed expression for the normalizing factor $B_k$ (see \cite{Er}, vol 1). $\qed$
\vsni
\noindent
Direct application of the Lemma \ref{pi1} now yields 
\begin{corollary} \label{P1}
The spectrum of ${\P}_1$ when acting upon ${\cal H}$ is given for each $r\in [0,1]$ by 
${\rm sp}\,({\P}_1)={\rm sp}\, (N_r)$. Each
eigenvalue ${\nu}_k$, with $k\geq 1$, is simple
and the corresponding (normalized) eigenfunction ${\xi}_k$ is given by 
\be
{\xi}_k (x)=({\B}_r\, [{\psi}_k])(x)=k\, B_k\,  {(1+(\beta_r-\alpha_r)x)^{k-1}\over
(1+\beta_r x)^{k+1}}.
\ee
\end{corollary}

\subsection{The spectrum of $\P$}
First, putting together (\ref{trM}), (\ref{trN}) and Theorem \ref{pi} we deduce the following result,
\begin{theorem} For all $r\in [0,1)$ the transfer operator $\P$ when acting upon ${\cal H}$ is of the trace-class, with
\be\label{tracciona}
{\rm tr}\, \P = {1\over \delta} + {\sqrt{1+4\rho} -1\over 2\,\sqrt{1+4\rho}}
\ee
with $\delta =1-r$ and $\rho =2-r$.
\end{theorem} 
\begin{remark}
Using (\ref{derfix1})  and setting  set $x_0\equiv 0$ one easily verifies that
\be
{1\over \delta} ={|\Phi'_{r,0}(x_0)|\over 1-\Phi'_{r,0}(x_0)}\, \cdot
\ee
Moreover, a straightforward computation based on the integral
\be
\int_0^\infty e^{-ax}\, J_1(bx) \, dx = {\sqrt{a^2+b^2} -a \over b\sqrt{a^2+b^2}}
\ee
and using (\ref{derfix1}) shows that 
\be
{\sqrt{1+4\rho} -1\over 2\,\sqrt{1+4\rho}} = {|\Phi'_{r,1}(x_1)|\over 1-\Phi'_{r,1}(x_1)}\, ,
\ee
where $x_1$, defined in (\ref{fixpt}), is the unique fixed point of $F_r$ besides $x_0$ (obviously one would arrive at the same conclusion
by directly summing the geometric series $ \sum_{k\geq 1} \mu_k + \sum_{k\geq 1}{\nu}_k$ with ${\mu}_k$ and ${\nu}_k$
given by (\ref{e1}) and (\ref{e2}), respectively).
As a result we can rewrite (\ref{tracciona})
as
\be
{\rm tr}\, \P = \sum_{i=0,1}{|\Phi'_{r,i}(x_i)|\over 1-\Phi'_{r,i}(x_i)}.
\ee
This expression does'nt come unexpectedly: it is an instance of a trace formula valid for more general analytic Markov maps 
(see \cite{Ma2}, Sec. 7.3.1).
\end{remark}

\vsni
\noindent
We are now ready to investigate the spectrum ${\rm sp}\, ({\P})$ on the space ${\cal H}$.
We first notice that from (\ref{formulona}) with $a=\sqrt{\rho}$ and $b=\sqrt{\rho}/2$ it follows that 
\be
N_r \left[ L_{k}^1 (2t)\right] =
{(-1)^k \over \rho} \, e^{-{r \over \rho}t}\,
  L_{k}^1\left({2t\over \rho}\right) = (-1)^k\, M_r \left[ L_{k}^1 (2t)\right] 
\ee
Therefore for all odd $k$ the function $L_{k}^1 (2t)$ lies in the kernel of $M_r+N_r$. This is in
agreement with (\ref{conj}) and (\ref{conj1}) since we have
\be
{\cal B}\left[ L_{k}^1 (2t)\right] =(k+1)(1-2x)^k
\ee
which for $k$ odd is an odd function w.r.t. $x=1/2$. This implies that $0$ is an eigenvalue
of infinite multiplicity for $\P$ for all $r\in [0,1]$.

\vsni
\noindent
Let us first consider the two extremal cases $r=0$ and $r=1$.
For $r=0$ we the above theorem gives ${\rm tr}\, \P =4/3$, suggesting the eigenvalues
$2^{-2n}$, ${n\geq 0}$. To check, we first note that
$$
M_0 \,t^{k} = {t^{k}\over 2^{k+1}}\quad\hbox{and}\quad
N_0 \,t^{k} = {k!\over 2}\,L_k^1(t/2)={(-t)^{k}\over 2^{k+1}}+ \sum_{j=0}^{k-1}d_j\, t^j
$$
with $d_j={(-1)^j\, (k+1)!\, k!\over 2^{j+1}\, (j+1)!\, j!\,  (k-j)!}\cdot$ Therefore the set of polynomials with even degree is a 
subset of $L^2(m)$ which is invariant for $M_0+N_0$, and any polynomial with odd degree is mapped into this subset. 
Given $k=2n$ it is now a simple task to construct a polynomial eigenfunction $\phi_{2n}$ (of degree $2n$) to the eigenvalue $2^{-2n}$.
The first three are 
$$
\phi_0(t)=1,\quad \phi_2(t)={t^2\over 2}-3t+2,\quad \phi_4(t)={t^4\over 24}-{5t^3\over 6}+{10t^2\over 3}-{32\over 15}\, \cdot
$$
Note moreover that for $\phi_{k}(t) =  \sum_{j=0}^{k}a_j\, t^j$ we have $({\B}\, [\phi_{k}])(x)= \sum_{j=0}^{k}a_j\,(j+1)!\, x^j$.
These observations along with standard arguments yield the following
\begin{proposition} 
For $r=0$ the spectrum of $\P$ when acting upon ${\cal H}$ is the set ${\rm sp}\, ({\P})=\{2^{-2n}\}_{n\geq 0}\cup\{0\}$. Each eigenvalue $2^{-2n}$ is simple
and the corresponding eigenfunction is a polynomial of degree $2n$.
\end{proposition}
At the opposite extremum we have the following result.
\begin{proposition} For $r=1$ the spectrum of $\P$ when acting upon ${\cal H}$ is 
 the union of
$[0,1]$ and a (possibily empty) countable set of real eigenvalues of finite multiplicity.
\end{proposition} 
{\sl Proof.} The assertion follows by noting that for $r=1$ the operator $\P$ 
when acting on ${\cal H}$ is isomorphic to $M_1+N_1$, 
which is a  self-adjoint compact perturbation of the (self-adjoint) multiplication operator $M_1$.
The assertion is now a consequence of Theorem 5.2 in \cite{GK}.
Note that although
the function $e_1(x)=(\log{2})^{-1}/x$ satisfies $\P e_1=e_1$ it does not belong to ${\cal H}$
and therefore $1\notin {\rm sp}\, ({\P}:{\cal H}\to {\cal H})$. $\qed$

\vsni
\noindent
We finally state the following
\begin{conjecture}
For all $0\leq r<1$ the spectrum of $\P$ when acting upon ${\cal H}$ is 
a countable subset of 
$[0,1]$ densely filling $[0,1]$ as $r\nearrow 1$. When $r=1$ the spectrum is purely continuous (i.e. there are no eigenvalues).
\end{conjecture}

\subsection{Operator-valued functions and zeta functions}
We may consider the operator-valued function ${\cal Q}_z$
defined as
\be\label{emmezeta}
{\cal Q}_z := z\, {\P_1}(1-z{\P_0})^{-1}.
\ee
Its relevance is twofold:
first, expanding formally in powers of $z$ we get
\be\label{qzeta}
{\cal Q}_z = \sum_{n=1}^\infty z^n \P_1 \P_0^{n-1}
\ee
and using (\ref{ind}) we see that for $z=1$ the operator ${\cal Q} \equiv {\cal Q}_1$ is the transfer operator
associated to the induced map $G_r$. Second,
it is related to ${\P}$ by the identity
\be
(1-{\cal Q}_z)(1-z\,\P_0) = 1-z\,\P.
\ee 
Now, as a consequence of (\ref{emmezeta}) and Lemmas \ref{pi0} - \ref{pi1}, we have the following expression
for the operator-valued function ${\cal Q}_z$ when acting on ${\cal H}$:
\be\label{Mop}
{\cal Q}_z {\B}\, [\varphi] =  {\B}\, \left[ N_r \left({1\over z}-\, M_r\right)^{-1}\varphi\right] \, \cdot
\ee
Putting together the above we obtain
\begin{theorem} For all $r\in [0,1]$
the operator-valued function ${{\cal Q}}_z$ when acting on the Hilbert space ${\cal H}$ is analytic
for $z\in \C \setminus \Lambda_r$, where $\Lambda_r = \{\rho^k\}_{k=1}^\infty$ for $r<1$ and $\Lambda_1=(1,\infty)$.
For each $z\in \C \setminus \Lambda_r$ it defines a trace-class operator.
\end{theorem}
\begin{remark}
Introducing the operator 
\be
L_r := (1-M_r)^{-1} \, N_r\, ,
\ee
we can rewrite (\ref{Pop}) and (\ref{Mop}) with $z=1$ as 
\be
{\P} {\B}\, [\varphi]  = {\B}\, \left[\, (M_r+(1-M_r)L_r\,) \,\varphi\right],
\ee
 and (recall that ${\cal Q} \equiv {\cal Q}_1$)
\be
{\cal Q} \,{\B}\, [\varphi] 
={\B}\, \left[(1-M_r) L_r (1-M_r)^{-1} \,\varphi\right],
\ee
respectively.
We thus see that the functions $\phi_r$ and $\psi_r$ defined in (\ref{densita}) satisfy
\be
L_r \phi_r  = \phi_r \quad \hbox{and}\quad
(1-M_r) L_r (1-M_r)^{-1}\psi_r =\psi_r,
\ee
so that 
\be
\phi_r = (1-M_r)^{-1} \psi_r\quad\hbox{and}\quad \psi_r=N_r \phi_r.
\ee

\end{remark}
We now consider the dynamical zeta functions $\zeta_{F_r}$ and $\zeta_{G_r}$
associated to the maps $F_r$ and $G_r$, respectively, and defined by 
the following formal series \cite{Ba2}:
\be
\zeta_{F_r} (z) = \exp \sum_{n=1}^{\infty} {z^n\over n} Z_n(F_r) \quad\hbox{and}\quad
\zeta_{G_r} (s) = \exp \sum_{n=1}^{\infty} {s^n\over n} Z_n(G_r),
\ee
where the `partition functions' $Z_n(F_r)$ and $Z_n(G_r)$ are given by
\be
Z_n(F_r) =\sum_{x=F_r^n(x)} \prod_{k=0}^{n-1}{1\over |F_r'(F_r^k(x))|}
\quad\hbox{and}\quad
Z_n(G_r) =\sum_{x=G_r^n(x)} \prod_{k=0}^{n-1}{1\over |G_r'(G_r^k(x))|}\cdot
\ee
Moreover, let us define the `grand partition function' 
\be\label{Xi}
\Xi_n(z) := 
\sum_{\ell=0}^{\infty}z^{\ell+n}\sum_{\scriptstyle x=G_r^n(x)=F_r^{\ell+n}(x)}
 \prod_{k=0}^{n-1}{1\over |G_r'(G_r^k(x))|},
\ee
and the
two-variable zeta function 
\be\label{twovar}
\z_2 (s,z) := \exp \sum_{n=1}^{\infty} {s^n\over n}\, \Xi_n(z).
\ee
A straightforward extension of (\cite{Is}, Proposition 4.3) to the present situation yields the identities
\be\label{due}
\z_2 (1,z) = (1- z)\, \z_{F_r} (z) \quad\hbox{and}\quad \z_2 (s,1) = \z_{G_r} (s)
\ee
which are valid for all $r\in [0,1]$ and wherever the series expansions converge absolutely.
Therefore the analytic properties of the dynamical zeta functions $\zeta_{F_r}$ and $\zeta_{G_r}$ can be deduced
from those of $\z_2 (s,z)$. In turn, the latter can be studied as follows.
For $q=0,1,\dots$ define
\begin{eqnarray}
{\P}_{0,q}f (x)&:=&\left[{\rho \over (\rho+rx)^2}\right]^{1+q}\cdot f \left({x \over \rho+rx}\right)\nonumber \\
{\P}_{1,q}f (x) &:=&\left[{\rho \over (\rho+rx)^2}\right]^{1+q}\cdot f \left(1-{x \over \rho+rx}\right)\nonumber
\end{eqnarray}
so that ${\P}_{0,0}\equiv {\P}_{0}$ and ${\P}_{1,0}\equiv {\P}_{1}$. These operators are supposed to act
upon the Hilbert space ${\cal H}_{q}\subseteq {\cal H}$ such that 
a function $f\in {\cal H}_{q}$ can be represented as
\be\label{represent3}
f (x)= ({\B}_q\, [\varphi])(x):={1\over x^{2(1+q)}}\int_0^\infty  e^{-{t\over x}}\,e^t\, 
\varphi (t)\, dm_q(t),\quad 
\varphi \in L^2(m_q),
\ee
with
\be
dm_q(t)= t^{2q+1}\,e^{-t}\, dt
\ee
A straightforward computation extending to non zero $q$ values those performed in the
previous Section yields 
\be \label{bundi}
{\P}_{0,q}\,{\B}_q\, [\varphi\, ] = {\B}_q\, [M_{r,q} \varphi\, ]\quad \hbox{and}\quad 
{\P}_{1,q}\,{\B}_q\, [\varphi\, ] = {\B}_q\, [N_{r,q} \varphi\, ]
\ee
where $M_{r,q} :L^2(m_q)\to L^2(m_q)$ is defined as
\be
(M_{r,q} \varphi )(t) := {1\over \rho^{1+q}}\,e^{-{r\over \rho}t}\, \varphi\left({t\over \rho}\right).
\ee 
and ${N}_{r,q} : L_2(m_q)\to L_2(m_q)$ is given by
\be\label{opn}
({N}_{r,q} \varphi) (t) :=  {1\over \rho^{1+q}}\,e^{{\delta\over \rho}t}\, \int_0^\infty 
{J_{2q+1}\left({2}\sqrt{{st/ \rho}}\right)\over (st/\rho)^{q+{1\over 2}}}\, \varphi (s) \, dm_q(s).
\ee
We can now generalize (\ref{qzeta}) defining
\be
{\cal Q}_{z,q} := \sum_{n=1}^\infty z^n \P_{1,q} \P_{0,q}^{n-1}.
\ee
In particular ${\cal Q}_{z,0}\equiv {\cal Q}_{z}$. Now set $\Lambda_{r,q} = \{\rho^{k+q}\}_{k=1}^\infty$ and note that
$\Lambda_{r,q} \subseteq \Lambda_{r}$ for all $q\geq 0$.
Reasoning as above and using (\ref{bundi}) we have that for any given $q=0,1\dots$ the operator valued function $z\to {\cal Q}_{z,q}$
when acting on ${\cal H}_{q}$  is analytic
for $z\in \C \setminus \Lambda_{r,q}$ and for each $z$ in this domain ${\cal Q}_{z,q}$ defines a trace-class operator with
\be
{\cal Q}_{z,q}\,{\B}_q\, [\varphi\, ] = {\B}_q\, [(-1)^q \,N_{r,q} \left({1\over z}-\, M_{r,q}\right)^{-1} \varphi\, ].
\ee
\vsni
\noindent
Furthermore, a straightforward adaptation of (\cite{Ma1}, Corollaries 4 and 5) to our 
$z$-dependent situation leads to 
the following expression for the grand partition function:
\be\label{traceformula}
\Xi_n(z) = {\rm tr}\,\, {\cal Q}_{z,0}^n - {\rm tr}\,\, {\cal Q}_{z,1}^n.
\ee
This trace formula along with standard arguments (see \cite{Ma1}) allow us to write the two-variables
zeta function (\ref{twovar}) as a ratio of Fredholm determinants,
\be
\z_2 (s,z) 
={{\rm det}\, (1-s\,{\cal Q}_{z,1}) 
\over {\rm det}\, (1-s\,{\cal Q}_{z,0}) },
\ee
where by definition 
\be
{\rm det}\,(1-s\,{\cal Q}_{z,q}) = \exp \left( -\sum_{n =1}^\infty {s^n\over n}\, {\rm tr}\,\, 
{\cal Q}^n_{z,q}
\right)
\ee
is in the sense of Grothendieck \cite{G}.
Putting together the above we obtain the following result from which the analytic properties of  
the dynamical zeta functions associated to the maps $F_r$ and $G_r$ can be readily deduced via (\ref{due}).

\begin{theorem}\label{twoanal}

\noindent
\begin{enumerate}
\item for each $s\in \C$, the function $\z_2 (s,z)$, considered as a
function of the variable $z$,
extends to a meromorphic function in $\{z\in \C : z\notin \Lambda_r\}$. Its poles are located among 
those $z$-values such
that ${\cal Q}_{z}:{\cal H} \to {\cal H}$  has $1/s$ as an eigenvalue;
\item for each $z\in \C \setminus \Lambda_r$, the function $\z_2 (s,z)$, considered as a
function of the variable $s$,
extends to a meromorphic function in  $\C$. Its poles are located among the inverses of the
eigenvalues of ${\cal Q}_{z}:{\cal H} \to {\cal H}$.
\end{enumerate}
\end{theorem}
\begin{remark}
The first statement with $s=1$ shows that for $r=1$ the function $\zeta_{F_r}(z)$ has a non-polar singularity
at $z=1$. This can be related to the non-analytic behaviour of the free energy at $\beta =0$ discussed in Section \ref{free}.
\end{remark}

\end{document}